\documentclass[12pt]{article}%
\usepackage{amsfonts}
\usepackage{sw20bams}
\usepackage{amsmath}
\usepackage{amssymb}
\usepackage{graphicx}%
\providecommand{\U}[1]{\protect\rule{.1in}{.1in}}
\begin{document}

\title{How Long Might We Wait at Random?}
\author{Steven Finch}
\date{August 27, 2022}
\maketitle

\begin{abstract}
In discrete time, customers arrive at random. \ Each waits until one of three
servers is available; each thereafter departs at random. \ We seek the
distribution of maximum line length of idle customers. \ Algebraic expressions
obtained for the two-server scenario do not appear feasible here. \ We also
review well-known distributional results for maximum wait time associated with
an M/M/$1$ queue and speculate about their generalization.

\end{abstract}

\footnotetext{Copyright \copyright \ 2019 \&\ 2022 by Steven R. Finch. All
rights reserved.}In queueing theory, maximum line length and maximum wait time
are different sides of the same coin, one spatial (Section 1) and the other
temporal (Section 2). \ 

\section{In Spatium}

Let $0<r<1$ and $0<p<3r$. \ Consider the Julia program:%
\[%
\begin{array}
[c]{l}%
\text{\texttt{u = 0}}\\
\text{\texttt{m = 0}}\\
\text{\texttt{for t=1:n}}\\
\text{\texttt{\ \ x = rand()%
$<$%
p \ \ \ \ \ \ \ \ \# x=1 means that an arrival occurs}}\\
\text{\texttt{\ \ y1 = rand()%
$<$%
r \ \ \ \ \ \ \ \# y1+y2+y3=3 means that three departures occur}}\\
\text{\texttt{\ \ y2 = rand()%
$<$%
r \ \ \ \ \ }}\\
\text{\texttt{\ \ y3 = rand()%
$<$%
r \ \ \ \ \ }}\\
\text{\texttt{\ \ if u==0 }}\\
\text{\texttt{\ \ \ \ u = x\ \ \ \ \ \ \ \ \ \ \ \ \ \ \ \ \ \ \ \ \ \ \ \ \#
increment is 1 or 0\ \ }}\\
\text{\texttt{\ \ else}}\\
\text{\texttt{\ \ \ \ if u==1}}\\
\text{\texttt{\ \ \ \ \ \ u = max(0,u+x-y1) \ \ \ \ \ \ \ \ \ \# increment is
1, 0 or -1}}\\
\text{\texttt{\ \ \ \ elseif u==2}}\\
\text{\texttt{\ \ \ \ \ \ u = max(0,u+x-y1-y2)\ \ \ \ \ \ \ \# increment is 1,
0, -1 or -2}}\\
\text{\texttt{\ \ \ \ else}}\\
\text{\texttt{\ \ \ \ \ \ u = max(0,u+x-y1-y2-y3)\ \ \ \ \# increment is 1, 0,
-1, -2 or -3}}\\
\text{\texttt{\ \ \ \ end}}\\
\text{\texttt{\ \ end}}\\
\text{\texttt{\ \ m = max(m,u)}}\\
\text{\texttt{end}}\\
\text{\texttt{return m}}%
\end{array}
\]
which simulates the maximum value of a Geo/Geo/$3$ queue with LAS-DA over $n$
time steps. \ The Boolean expressions containing Julia's Uniform $[0,1]$
random deviate generator ensure that $X\sim\,$Bernoulli($p$) and $Y\sim
\,$Bernoulli($r$). \ The word \textquotedblleft Geometric\textquotedblright%
\ arises because%
\[%
\begin{array}
[c]{ccc}%
\mathbb{P}\left\{  \text{time lapse between adjacent arrivals is }i\right\}
=p\,q^{i-1}, &  & i\geq1
\end{array}
\]
where $q=1-p$ and
\[%
\begin{array}
[c]{ccc}%
\mathbb{P}\left\{  \text{time lapse between adjacent departures is }j\right\}
=r\,s^{j-1}, &  & j\geq1
\end{array}
\]
where $s=1-r$. \ Clearly $3s-2<q<1$. \ LAS\ stands for \textquotedblleft late
arrival system\textquotedblright\ and DA\ stands for \textquotedblleft delayed
access\textquotedblright\ \cite{Htr-heu}; in particular, a customer entering
an empty queue at time $t$ is not immediately eligible for service, but rather
at time $t+1$. We study the asymptotic distribution of the maximum $M_{n}$.
\ This section is best thought of as an extended addendum to \cite{Fnh-heu},
written for the sake of completeness. \ The only difference with our earlier
work is that closed-form expressions here become unwieldy and thus our
approach is more numeric and less symbolic.

The Poisson clumping heuristic \cite{Ald-heu}, while not a theorem, gives
results identical to exact asymptotic expressions when such exist, and
evidently provides excellent predictions otherwise. \ Consider an irreducible
positive recurrent Markov chain with stationary distribution $\pi$. \ For
large enough $k$, the maximum of the chain satisfies%
\[
\mathbb{P}\left\{  M_{n}<k\right\}  \sim\exp\left(  -\frac{\pi_{k}}%
{\mathbb{E}(C)}n\right)
\]
as $n\rightarrow\infty$, where $C$ is the sojourn time in $k$ during a clump
of nearby visits to $k$.

Starting with transition matrix%
\[
\left(  {\small
\begin{array}
[c]{cccccccc}%
q & p & 0 & 0 & 0 & 0 & 0 & \cdots\\
qr & pr+qs & ps & 0 & 0 & 0 & 0 & \cdots\\
qr^{2} & pr^{2}+2qrs & 2prs+qs^{2} & ps^{2} & 0 & 0 & 0 & \cdots\\
qr^{3} & pr^{3}+3qr^{2}s & 3pr^{2}s+3qrs^{2} & 3prs^{2}+qs^{3} & ps^{3} & 0 &
0 & \cdots\\
0 & qr^{3} & pr^{3}+3qr^{2}s & 3pr^{2}s+3qrs^{2} & 3prs^{2}+qs^{3} & ps^{3} &
0 & \cdots\\
0 & 0 & qr^{3} & pr^{3}+3qr^{2}s & 3pr^{2}s+3qrs^{2} & 3prs^{2}+qs^{3} &
ps^{3} & \cdots\\
0 & 0 & 0 & qr^{3} & pr^{3}+3qr^{2}s & 3pr^{2}s+3qrs^{2} & 3prs^{2}+qs^{3} &
\cdots\\
\vdots & \vdots & \vdots & \vdots & \vdots & \vdots & \vdots & \ddots
\end{array}
}\right)
\]
we obtain \cite{AH-heu}%
\[%
\begin{array}
[c]{ccccc}%
\pi_{j}=\pi_{j}^{\ast}\pi_{3} &  & \text{for} &  & 0\leq j\leq2,
\end{array}
\]%
\[
\pi_{3}=\left(  \frac{1}{1-\omega}+\pi_{0}^{\ast}+\pi_{1}^{\ast}+\pi_{2}%
^{\ast}\right)  ^{-1},
\]%
\[%
\begin{array}
[c]{ccccc}%
\pi_{j}=\omega^{j-3}\pi_{3} &  & \text{for} &  & j\geq4
\end{array}
\]
where%
\[
\pi_{2}^{\ast}=\frac{1}{ps^{2}}\left[  \left(  1-3prs^{2}-qs^{3}\right)
-\left(  3pr^{2}s+3qrs^{2}\right)  \omega-\left(  pr^{3}+3qr^{2}s\right)
\omega^{2}-\left(  qr^{3}\right)  \omega^{3}\right]  ,
\]%
\[
\pi_{1}^{\ast}=\frac{1}{ps}\left[  \left(  1-2prs-qs^{2}\right)  \pi_{2}%
^{\ast}-\left(  3pr^{2}s+3qrs^{2}\right)  -\left(  pr^{3}+3qr^{2}s\right)
\omega-\left(  qr^{3}\right)  \omega^{2}\right]  ,
\]%
\[
\pi_{0}^{\ast}=\frac{1}{p}\left[  \left(  1-pr-qs\right)  \pi_{1}^{\ast
}-\left(  pr^{2}+2qrs\right)  \pi_{2}^{\ast}-\left(  pr^{3}+3qr^{2}s\right)
-\left(  qr^{3}\right)  \omega\right]
\]
and $0<\omega<1$ satisfies the cubic equation $\omega=(q\omega+p)(r\omega
+s)^{3}$, that is,%
\[
\omega=\frac{1}{3qr}\left[  \left(  -3+2r+3ps\right)  +(3q-r)r\left(  \frac
{2}{\chi+\theta}\right)  ^{1/3}-\left(  \frac{\chi+\theta}{2}\right)
^{1/3}\right]
\]
where%
\[%
\begin{array}
[c]{ccc}%
\chi=-9qr^{2}+2r^{3}-27q^{2}s, &  & \theta=\sqrt{4(3q-r)^{3}r^{3}+\chi^{2}}.
\end{array}
\]
Note that, if $k=\log_{1/\omega}(n)+h+1$, we have%
\[
\left(  \dfrac{1}{\omega}\right)  ^{k}=n\left(  \dfrac{1}{\omega}\right)
^{h+1}%
\]
thus%
\[
\pi_{k}n=\pi_{3}\omega^{k-3}n=\pi_{3}\omega^{h-2}.
\]

We need now to calculate $\mathbb{E}(C)$. Consider a random walk on the
integers consisting of incremental steps satisfying%
\[
\left\{
\begin{array}
[c]{lll}%
-3 &  & \text{with probability }qr^{3},\\
-2 &  & \text{with probability }pr^{3}+3qr^{2}s,\\
-1 &  & \text{with probability }3pr^{2}s+3qrs^{2},\\
0 &  & \text{with probability }3prs^{2}+qs^{3},\\
1 &  & \text{with probability }ps^{3}.
\end{array}
\right.
\]
For nonzero $j$, let $\nu_{j}$ denote the probability that, starting from
$-j$, the walker eventually hits $0$. \ Let $\nu_{0}$ denote the probability
that, starting from $0$, the walker eventually returns to $0$ (at some future
time). \ We have two values for $\nu_{0}$: when it is used in a recursion, it
is equal to $1$; when it corresponds to a return probability, it retains the
symbol $\nu_{0}$. \ Let $j\geq1$. \ Using%
\[
\nu_{j}=ps^{3}\nu_{j-1}+(3prs^{2}+qs^{3})\nu_{j}+(3pr^{2}s+3qrs^{2})\nu
_{j+1}+(pr^{3}+3qr^{2}s)\nu_{j+2}+qr^{3}\nu_{j+3},
\]%
\[
\nu_{0}=ps^{3}\nu_{-1}+(3prs^{2}+qs^{3})+(3pr^{2}s+3qrs^{2})\nu_{1}%
+(pr^{3}+3qr^{2}s)\nu_{2}+qr^{3}\nu_{3}%
\]
define%
\begin{align*}
F(z)  &  =%
{\displaystyle\sum\limits_{j=1}^{\infty}}
\nu_{j}z^{j}\\
&  =ps^{3}z%
{\displaystyle\sum\limits_{j=1}^{\infty}}
\nu_{j-1}z^{j-1}+(3prs^{2}+qs^{3})%
{\displaystyle\sum\limits_{j=1}^{\infty}}
\nu_{j}z^{j}+\frac{3pr^{2}s+3qrs^{2}}{z}%
{\displaystyle\sum\limits_{j=1}^{\infty}}
\nu_{j+1}z^{j+1}\\
&  +\frac{pr^{3}+3qr^{2}s}{z^{2}}%
{\displaystyle\sum\limits_{j=1}^{\infty}}
\nu_{j+2}z^{j+2}+\frac{qr^{3}}{z^{3}}%
{\displaystyle\sum\limits_{j=1}^{\infty}}
\nu_{j+3}z^{j+3}\\
&  =ps^{3}z\left[  F(z)+1\right]  +(3prs^{2}+qs^{3})F(z)+\frac{3pr^{2}%
s+3qrs^{2}}{z}\left[  F(z)-\nu_{1}z\right] \\
&  +\frac{pr^{3}+3qr^{2}s}{z^{2}}\left[  F(z)-\nu_{1}z-\nu_{2}z^{2}\right]
+\frac{qr^{3}}{z^{3}}\left[  F(z)-\nu_{1}z-\nu_{2}z^{2}-\nu_{3}z^{3}\right]
\end{align*}
equivalently%
\begin{align*}
&  \left[  1-ps^{3}z-3prs^{2}-qs^{3}-\frac{3pr^{2}s+3qrs^{2}}{z}-\frac
{pr^{3}+3qr^{2}s}{z^{2}}-\frac{qr^{3}}{z^{3}}\right]  F(z)\\
&  =ps^{3}z-\frac{3pr^{2}s+3qrs^{2}}{z}(\nu_{1}z)-\frac{pr^{3}+3qr^{2}s}%
{z^{2}}(\nu_{1}z+\nu_{2}z^{2})-\frac{qr^{3}}{z^{3}}(\nu_{1}z+\nu_{2}z^{2}%
+\nu_{3}z^{3})
\end{align*}
equivalently%
\begin{align*}
&  \left[  qr^{3}+(pr^{3}+3qr^{2}s)z+(3pr^{2}s+3qrs^{2})z^{2}-(1-3prs^{2}%
-qs^{3})z^{3}+ps^{3}z^{4}\right]  F(z)\\
&  =-ps^{3}z^{4}+(3pr^{2}s+3qrs^{2})z^{2}(\nu_{1}z)+(pr^{3}+3qr^{2}s)z(\nu
_{1}z+\nu_{2}z^{2})+qr^{3}(\nu_{1}z+\nu_{2}z^{2}+\nu_{3}z^{3})
\end{align*}
equivalently%
\begin{align*}
&  (1-z)\left[  qr^{3}+(3q-2r+3pr)r^{2}z+(1+s+s^{2}-3ps^{2})rz^{2}-ps^{3}%
z^{3}\right]  F(z)\\
&  =-ps^{3}z^{4}+3pr^{2}sz^{3}\nu_{1}+3qrs^{2}z^{3}\nu_{1}+pr^{3}z^{2}\nu
_{1}+pr^{3}z^{3}\nu_{2}+3qr^{2}sz^{2}\nu_{1}+3qr^{2}sz^{3}\nu_{2}+qr^{3}%
zv_{1}\\
&  +qr^{3}z^{2}\nu_{2}+z^{3}\left[  \nu_{0}-ps^{3}\nu_{-1}-(3prs^{2}%
+qs^{3})-(3pr^{2}s+3qrs^{2})\nu_{1}-(pr^{3}+3qr^{2}s)\nu_{2}\right] \\
&  =z^{3}\nu_{0}+\left(  qr+prz+3qsz\right)  r^{2}z\nu_{1}+qr^{3}z^{2}\nu
_{2}-ps^{3}z^{3}\nu_{-1}-3prs^{2}z^{3}-qs^{3}z^{3}-ps^{3}z^{4}.
\end{align*}
Examine the denominator of $F(z)$. \ Only the first three of its four zeroes
$z_{1}$, $z_{2}$, $1$, $z_{3}$ are of interest (the fourth is $>1$).
\ Substituting $z=1$, $z=z_{1}$ and $z=z_{2}$ into the numerator $N_{F}$ of
$F(z)$, then setting $N_{F}=0$, gives three equations in four unknowns.\ 

Let $j\geq1$. \ Using%
\[
\nu_{-j}=ps^{3}\nu_{-j-1}+(3prs^{2}+qs^{3})\nu_{-j}+(3pr^{2}s+3qrs^{2}%
)\nu_{-j+1}+(pr^{3}+3qr^{2}s)\nu_{-j+2}+qr^{3}\nu_{-j+3},
\]%
\[
\nu_{0}=ps^{3}\nu_{-1}+(3prs^{2}+qs^{3})+(3pr^{2}s+3qrs^{2})\nu_{1}%
+(pr^{3}+3qr^{2}s)\nu_{2}+qr^{3}\nu_{3}%
\]
define%
\begin{align*}
G(z)  &  =%
{\displaystyle\sum\limits_{j=1}^{\infty}}
\nu_{-j}z^{j}\\
&  =\frac{ps^{3}}{z}%
{\displaystyle\sum\limits_{j=1}^{\infty}}
\nu_{-j-1}z^{j+1}+(3prs^{2}+qs^{3})%
{\displaystyle\sum\limits_{j=1}^{\infty}}
\nu_{-j}z^{j}+(3pr^{2}s+3qrs^{2})z%
{\displaystyle\sum\limits_{j=1}^{\infty}}
\nu_{-j+1}z^{j-1}\\
&  +(pr^{3}+3qr^{2}s)z^{2}%
{\displaystyle\sum\limits_{j=1}^{\infty}}
\nu_{-j+2}z^{j-2}+qr^{3}z^{3}%
{\displaystyle\sum\limits_{j=1}^{\infty}}
\nu_{-j+3}z^{j-3}\\
&  =\frac{ps^{3}}{z}\left[  G(z)-\nu_{-1}z\right]  +(3prs^{2}+qs^{3}%
)G(z)+(3pr^{2}s+3qrs^{2})z\left[  G(z)+1\right] \\
&  +(pr^{3}+3qr^{2}s)z^{2}\left[  G(z)+1+\frac{\nu_{1}}{z}\right]
+qr^{3}z^{3}\left[  G(z)+1+\frac{\nu_{1}}{z}+\frac{\nu_{2}}{z^{2}}\right]
\end{align*}
equivalently%
\begin{align*}
&  \left[  1-\frac{ps^{3}}{z}-3prs^{2}-qs^{3}-(3pr^{2}s+3qrs^{2}%
)z-(pr^{3}+3qr^{2}s)z^{2}-qr^{3}z^{3}\right]  G(z)\\
&  =-\frac{ps^{3}}{z}(\nu_{-1}z)+(3pr^{2}s+3qrs^{2})z+(pr^{3}+3qr^{2}%
s)z^{2}\left(  1+\frac{\nu_{1}}{z}\right)  +qr^{3}z^{3}\left(  1+\frac{\nu
_{1}}{z}+\frac{\nu_{2}}{z^{2}}\right)
\end{align*}
equivalently%
\begin{align*}
&  \left[  qr^{3}z^{4}+(pr^{3}+3qr^{2}s)z^{3}+(3pr^{2}s+3qrs^{2}%
)z^{2}-(1-3prs^{2}-qs^{3})z+ps^{3}\right]  G(z)\\
&  =ps^{3}(\nu_{-1}z)-(3pr^{2}s+3qrs^{2})z^{2}-(pr^{3}+3qr^{2}s)z^{2}%
(z+\nu_{1})-qr^{3}z^{2}\left(  z^{2}+\nu_{1}z+\nu_{2}\right)
\end{align*}
equivalently%
\begin{align*}
&  (1-z)\left[  ps^{3}-(1+s+s^{2}-3ps^{2})rz-(3q-2r+3pr)r^{2}z^{2}-qr^{3}%
z^{3}\right]  G(z)\\
&  =ps^{3}z\nu_{-1}-(3prs+3qs^{2}+pr^{2}z+3qrsz+qr^{2}z^{2})rz^{2}%
-(pr+3qs+qrz)r^{2}z^{2}\nu_{1}-qr^{3}z^{2}\nu_{2}.
\end{align*}
Examine the denominator of $G(z)$. \ Only the zero $z_{4}$ of smallest modulus
interests us. \ Substituting $z=z_{4}$ into the numerator $N_{G}$ of $G(z)$,
and setting $N_{G}=0$, gives a fourth equation (to include with the other
three from earlier). \ For example, if%
\[%
\begin{array}
[c]{ccccccc}%
p=\frac{1}{3}, &  & q=\frac{2}{3}, &  & r=\frac{1}{6}, &  & s=\frac{5}{6}%
\end{array}
\]
we have\
\[%
\begin{array}
[c]{ccc}%
\omega=0.5744080010..., &  & \pi_{3}=0.1777380492...
\end{array}
\]
from earlier and%
\[%
\begin{array}
[c]{ccc}%
\nu_{0}=0.8437587438..., &  & \nu_{-1}=0.9309681530...,\\
\nu_{1}=0.5744080010..., &  & \nu_{2}=0.3299445517...
\end{array}
\]
after solving the simultaneous system in $\nu_{0}$, $\nu_{-1}$, $\nu_{1}$,
$\nu_{2}$. \ Observe that%
\begin{align*}
\mathbb{P}\left\{  M_{n}\leq\log_{1/\omega}(n)+h\right\}   &  =P\left\{
M_{n}<\log_{1/\omega}(n)+h+1\right\} \\
&  \sim\exp\left[  -\dfrac{\pi_{3}(1-\nu_{0})}{\omega^{2}}\omega^{h}\right]
\end{align*}
as $n\rightarrow\infty$. \ The ratio within the exponential argument is%
\[
\dfrac{\pi_{3}(1-\nu_{0})}{\omega^{2}}=0.0841657058...
\]
and hence%
\begin{align*}
\mathbb{E}\left(  M_{n}\right)   &  \approx\frac{\ln(n)}{\ln(\frac{1}{\omega
})}+\frac{\gamma+\ln\left(  \tfrac{\pi_{3}(1-\nu_{0})}{\omega^{2}}\right)
}{\ln(\frac{1}{\omega})}+\frac{1}{2}\\
&  \approx(1.8037019224...)\ln(n)-(2.9229790566...)
\end{align*}
for sufficiently large $n$, where $\gamma$ denotes Euler's constant
\cite{Fi1-heu}. \ Such moment formulas usually contain tiny periodic
fluctuations, but we omit these from consideration.

The use of an expected maximum for performance analysis, instead of a simple
average, does not appear to lead to surprising outcomes. \ A\ corollary of the
preceding numerical results is that, in a busy hospital emergency room (with
$p=1/3$), one fast doctor (with $r=1/2$) outperforms three slow doctors (each
with $r=1/6$). \ For average queue lengths \cite{AH-heu},
\[%
{\displaystyle\sum\limits_{j=1}^{\infty}}
j\pi_{j}=\frac{pq}{r-p}=1.33333...
\]
corresponding to Geo/Geo/$1$ and
\[%
{\displaystyle\sum\limits_{j=1}^{\infty}}
j\pi_{j}=\pi_{1}+2\pi_{2}+\frac{3-2\omega}{(1-\omega)^{2}}\pi_{3}=2.56365...
\]
corresponding to Geo/Geo/$3$. \ This is also consistent with results in
\cite{Fi2-heu} governing deterministic traffic signals: we do better with an
$RGRGRG...$ pattern than with $RRRGGG...$.

\section{In Tempore}

Let $c\geq1$ be an integer and $0<\lambda<c\mu$. \ Consider the R program:%
\[%
\begin{array}
[c]{l}%
\text{\texttt{K
$<$%
- rpois(1,n*lambda)}}\\
\text{\texttt{P
$<$%
- matrix(0,K,3) \ \ \ \ \ \ \ \ \# matrix of patients}}\\
\text{\texttt{P[,1]
$<$%
- sort(runif(K,0,n)) \ \ \ \ \ \ \ \ \# arrival times}}\\
\text{\texttt{P[,3]
$<$%
- rexp(K,mu) \ \ \ \ \ \ \ \ \ \ \ \ \ \ \ \ \# treatment lengths}}\\
\text{\texttt{D
$<$%
- rep(0,c) \ \ \ \ \ \ \ \ \ \ \ \ \ \# vector of doctors}}\\
\text{\texttt{k.sys
$<$%
- function(i,P) P[i,2]+P[i,3]-P[i,1]}}\\
\text{\texttt{k.que
$<$%
- function(i,P) P[i,2]-P[i,1]}}\\
\text{\texttt{for (i in 1:K)}}\\
\text{\texttt{\ \ \{}}\\
\text{\texttt{\ \ \ \ j
$<$%
- which.min(D)}}\\
\text{\texttt{\ \ \ \ P[i,2]
$<$%
- max(P[i,1],D[j])}}\\
\text{\texttt{\ \ \ \ D[j]
$<$%
- P[i,2] + P[i,3] \ \ \ \ \ \ \ \ \# departure times}}\\
\text{\texttt{\ \ \}}}\\
\text{\texttt{W.sys
$<$%
- sapply(1:K,k.sys,P=P)}}\\
\text{\texttt{W.que
$<$%
- sapply(1:K,k.que,P=P)}}\\
\text{\texttt{list(max(W.sys),max(W.que))}}%
\end{array}
\]
which simulates the maximum wait time associated with an M/M/$c$ queue over
the time interval $[0,n]$. \ More precisely, at any arrival time $t=P_{i,1}$,
let $W_{sys}$ denote the $i^{\text{th}}$ patient wait time in the
\textbf{system} (either queue or treatment) and $W_{que}$ denote the
$i^{\text{th}}$ patient wait time in the \textbf{queue} (excluding treatment).
\ For $c=1$, the maximums of $W_{sys}$ and $W_{que}$ over all arrival times up
to $n$ satisfy \cite{Ig1-heu, Ig2-heu, BW-heu}%
\[
\lim_{n\rightarrow\infty}\mathbb{P}\left\{  (\mu-\lambda)\max\limits_{0\leq
t\leq n}W_{sys}\leq\ln\left[  \lambda(1-\rho)^{2}n\right]  +x\right\}
=\exp(-\exp(-x)),
\]%
\[
\lim_{n\rightarrow\infty}\mathbb{P}\left\{  (\mu-\lambda)\max\limits_{0\leq
t\leq n}W_{que}\leq\ln\left[  \lambda(1-\rho)^{2}\rho\,n\right]  +x\right\}
=\exp(-\exp(-x))
\]
where $\rho=\lambda/\mu$. \ These imply%
\[
\mathbb{P}\left\{  \max\limits_{0\leq t\leq n}W_{sys}\leq y\right\}
\approx\exp\left[  -\lambda(1-\rho)^{2}n\exp\left[  -(\mu-\lambda)\right]
y\right]  ,
\]%
\[
\mathbb{P}\left\{  \max\limits_{0\leq t\leq n}W_{que}\leq y\right\}
\approx\exp\left[  -\lambda(1-\rho)^{2}\rho\,n\exp\left[  -(\mu-\lambda
)\right]  y\right]
\]
for large enough $n$. \ Also, if $\lambda=1/3$ and $\mu=1/2$, then%
\[
\mathbb{E}\left(  \max\limits_{0\leq t\leq n}W_{sys}\right)  =\frac{\ln
(n)}{\mu-\lambda}+\frac{\gamma+\ln\left[  \lambda(1-\rho)^{2}\right]  }%
{\mu-\lambda}=6\ln(n)-16.31172...,
\]%
\[
\mathbb{E}\left(  \max\limits_{0\leq t\leq n}W_{que}\right)  =\frac{\ln
(n)}{\mu-\lambda}+\frac{\gamma+\ln\left[  \lambda(1-\rho)^{2}\rho\right]
}{\mu-\lambda}=6\ln(n)-18.74451....
\]
In particular, assuming $n=20000$, the expected maximum wait times are
$43.109$ and $40.676$ respectively.

For $c\geq2$, no analogous formulas are known. \ Assuming $n=20000$,
$\lambda=1/3$ and $\mu=1/(2c)$, we estimate via simulation that%
\[%
\begin{array}
[c]{ccc}%
\mathbb{E}\left(  \max\limits_{0\leq t\leq n}W_{sys}\right)  \approx
51.0>43.109, &  & \mathbb{E}\left(  \max\limits_{0\leq t\leq n}W_{que}\right)
\approx39.4<40.676
\end{array}
\]
for $c=2$ and%
\[%
\begin{array}
[c]{ccc}%
\mathbb{E}\left(  \max\limits_{0\leq t\leq n}W_{sys}\right)  \approx
64.1>51.0, &  & \mathbb{E}\left(  \max\limits_{0\leq t\leq n}W_{que}\right)
\approx38.3<39.4
\end{array}
\]
for $c=3$. \ With regard to the expected maximum of $W_{sys}$, one fast doctor
outperforms $c$ slow doctors. \ With regard to the expected maximum of
$W_{que}$, however, $c$ slow doctors outperform one fast doctor.
\ Experimental histograms of such maximums appear to be roughly Gumbel-shaped:
a two-moment fit provides a fairly compelling but ultimately inconclusive
match between theory and data. \ Only rigorous analysis for $c\geq2$ will
clarify the situation.

The same $W_{que}$ trend occurs for simple averages. \ Given $\lambda=1/3$,
$\mu=1/(2c)$ as before and allowing $n\rightarrow\infty$, we have the
following when $c=1$:
\[
\mathbb{E}\left(  W_{que}\right)  =\mathbb{E}\left(  W_{sys}\right)  -\frac
{1}{\mu}=\frac{\lambda}{\left(  \mu-\lambda\right)  \mu}=4,
\]
when $c=2$:
\[
\mathbb{E}\left(  W_{que}\right)  =\mathbb{E}\left(  W_{sys}\right)  -\frac
{1}{\mu}=\frac{\lambda^{2}}{\left(  2\mu-\lambda\right)  \left(  2\mu
+\lambda\right)  \mu}=3.2
\]
and when $c=3$:
\[
\mathbb{E}\left(  W_{que}\right)  =\mathbb{E}\left(  W_{sys}\right)  -\frac
{1}{\mu}=\frac{\lambda^{3}}{\left(  3\mu-\lambda\right)  \left(  \lambda
^{2}+4\lambda\mu+6\mu^{2}\right)  \mu}=2.66666....
\]
Such mean results apply not only to the FIFO\ discipline (first in, first
out), but also to the LIFO\ (last in, first out) and SIRO (serve in random
order)\ disciplines. \ The underlying SIRO\ wait time distribution is more
complicated than that for FIFO, as shown for $c=1$ in \cite{Flt-heu, Lbv-heu}.

Alternative treatments of M/M/$1$ maximum wait times include \cite{Gm-heu,
AA-heu}, where $n$ or $k$ refer not to the length of the time interval but
instead the total number of customers. \ Simple averages of wait times for
Geo/Geo/$c$ can be calculated by methods given in \cite{Alf-heu}; maximums
remain open.

\section{Acknowledgements}

I am thankful to Guy Louchard for introducing me to the Poisson clumping
heuristic (especially recursions for $\nu_{j}$ and $\nu_{-j}$), and to Stephan
Wagner for extracting discrete Gumbel asymptotics in \cite{Fi3-heu} (a
contribution leading to \cite{Fnh-heu, Fi0-heu, Fic-heu}). \ See also
\cite{Fi4-heu, Fi5-heu} for work involving the case of Uniform$[a,b]$ service
times and Deterministic$[a]$ subcase ($b\rightarrow a$). \ The creators of R,
Julia, Mathematica and Matlab, as well as administrators of the MIT\ Engaging
Cluster and the MIT\ Supercloud Cluster, earn my gratitude every day.%
\begin{figure}
[ptb]
\begin{center}
\includegraphics[
height=4.7513in,
width=6.0122in
]%
{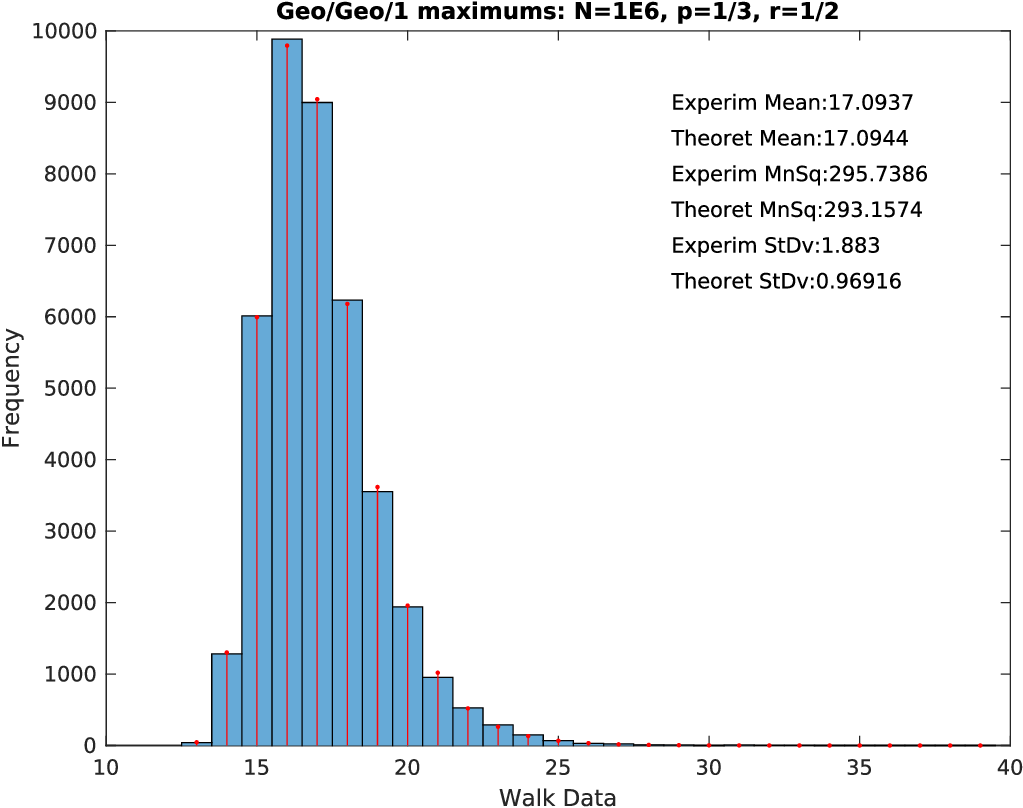}%
\caption{Maximum line length over $n$ time steps, with $1$ server}%
\end{center}
\end{figure}
\begin{figure}
[ptb]
\begin{center}
\includegraphics[
height=4.7859in,
width=6.0537in
]%
{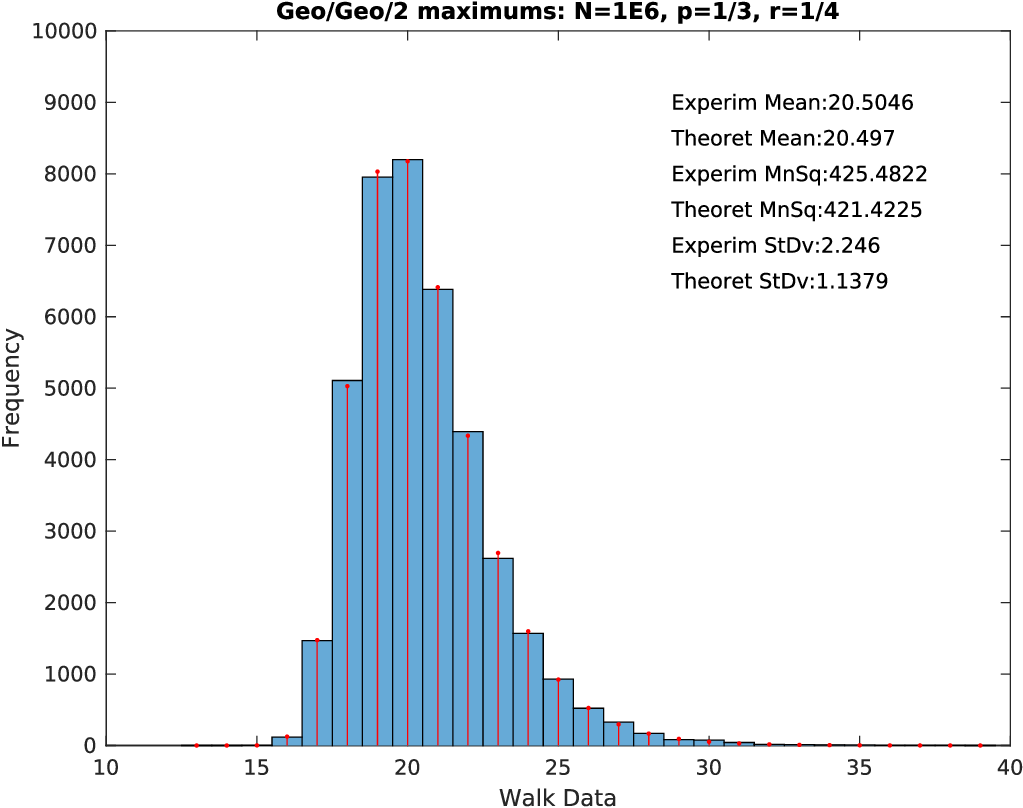}%
\caption{Maximum line length over $n$ time steps, with $2$ servers}%
\end{center}
\end{figure}
\begin{figure}
[ptb]
\begin{center}
\includegraphics[
height=4.7859in,
width=6.0537in
]%
{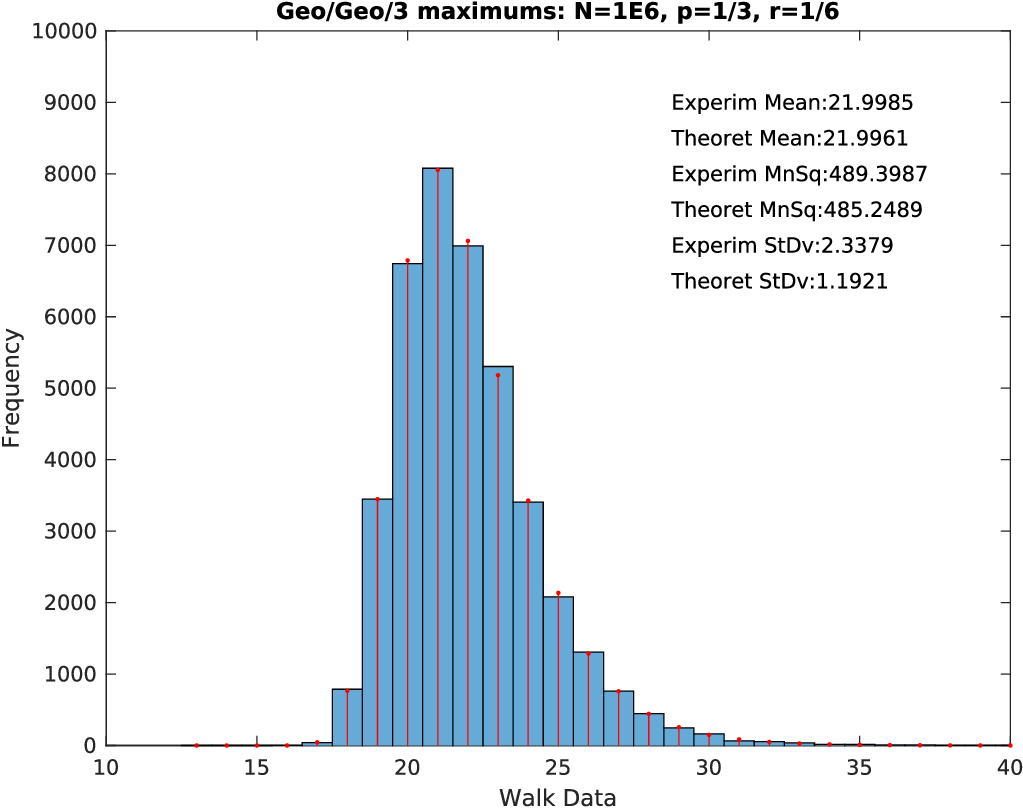}%
\caption{Maximum line length over $n$ time steps, with $3$ servers}%
\end{center}
\end{figure}
\begin{figure}
[ptb]
\begin{center}
\includegraphics[
height=3.0234in,
width=3.4541in
]%
{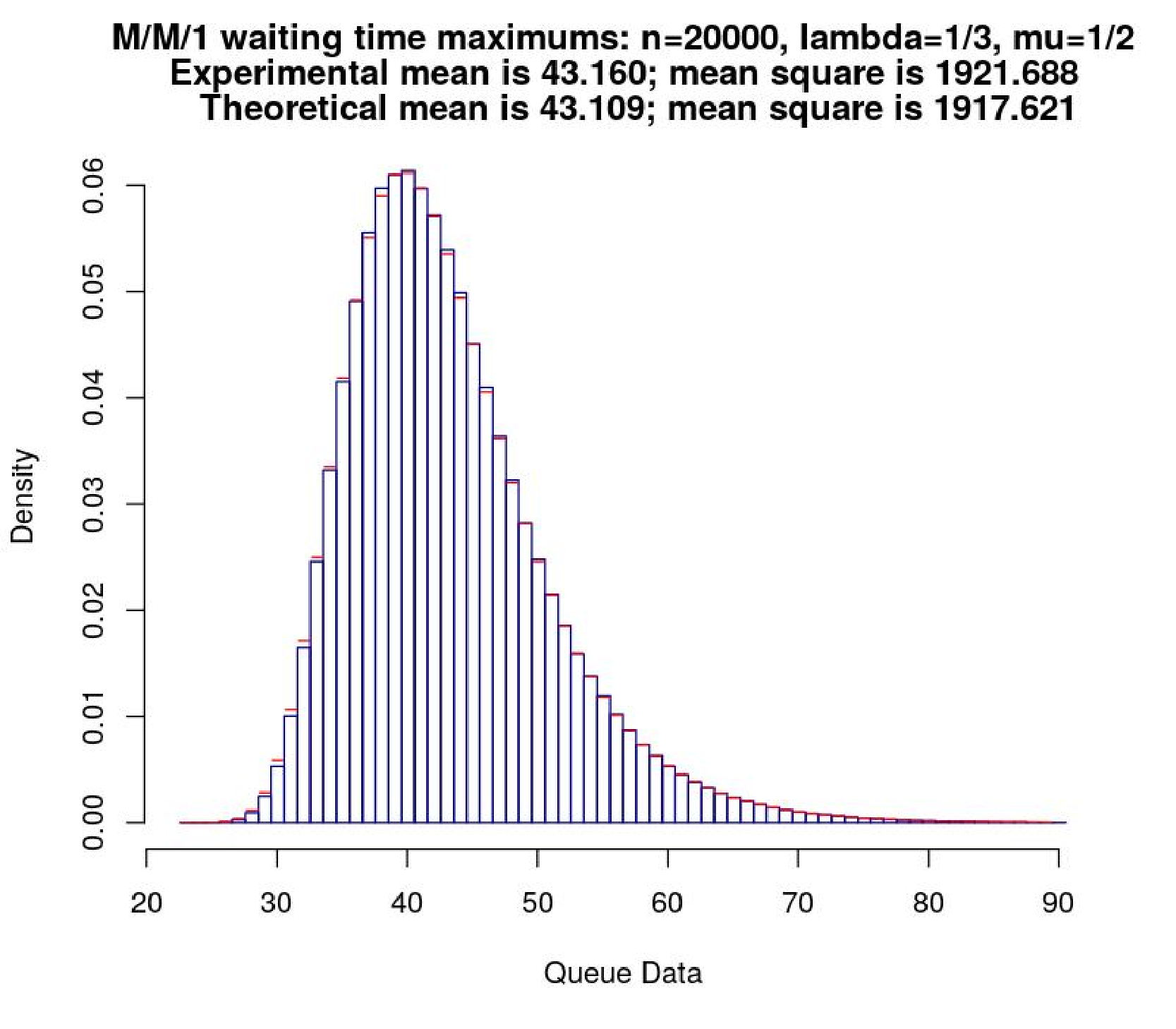}%
\caption{Maximum of $W_{sys}$ over time interval of length $n$, with 1 server}%
\end{center}
\end{figure}
\begin{figure}
[ptb]
\begin{center}
\includegraphics[
height=3.0234in,
width=3.4541in
]%
{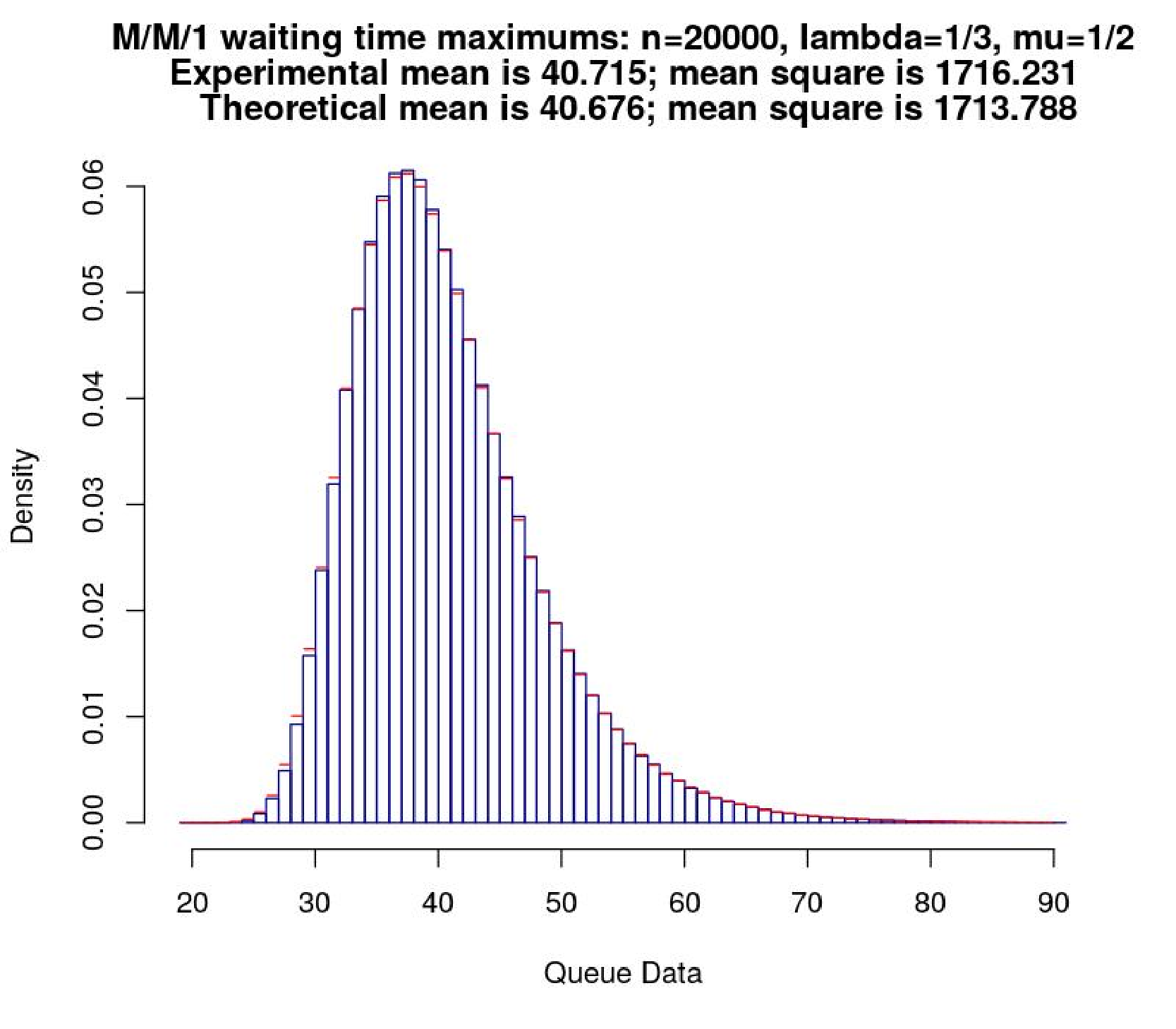}%
\caption{Maximum of $W_{que}$ over time interval of length $n$, with $1$
server}%
\end{center}
\end{figure}
\begin{figure}
[ptb]
\begin{center}
\includegraphics[
height=3.0234in,
width=3.4541in
]%
{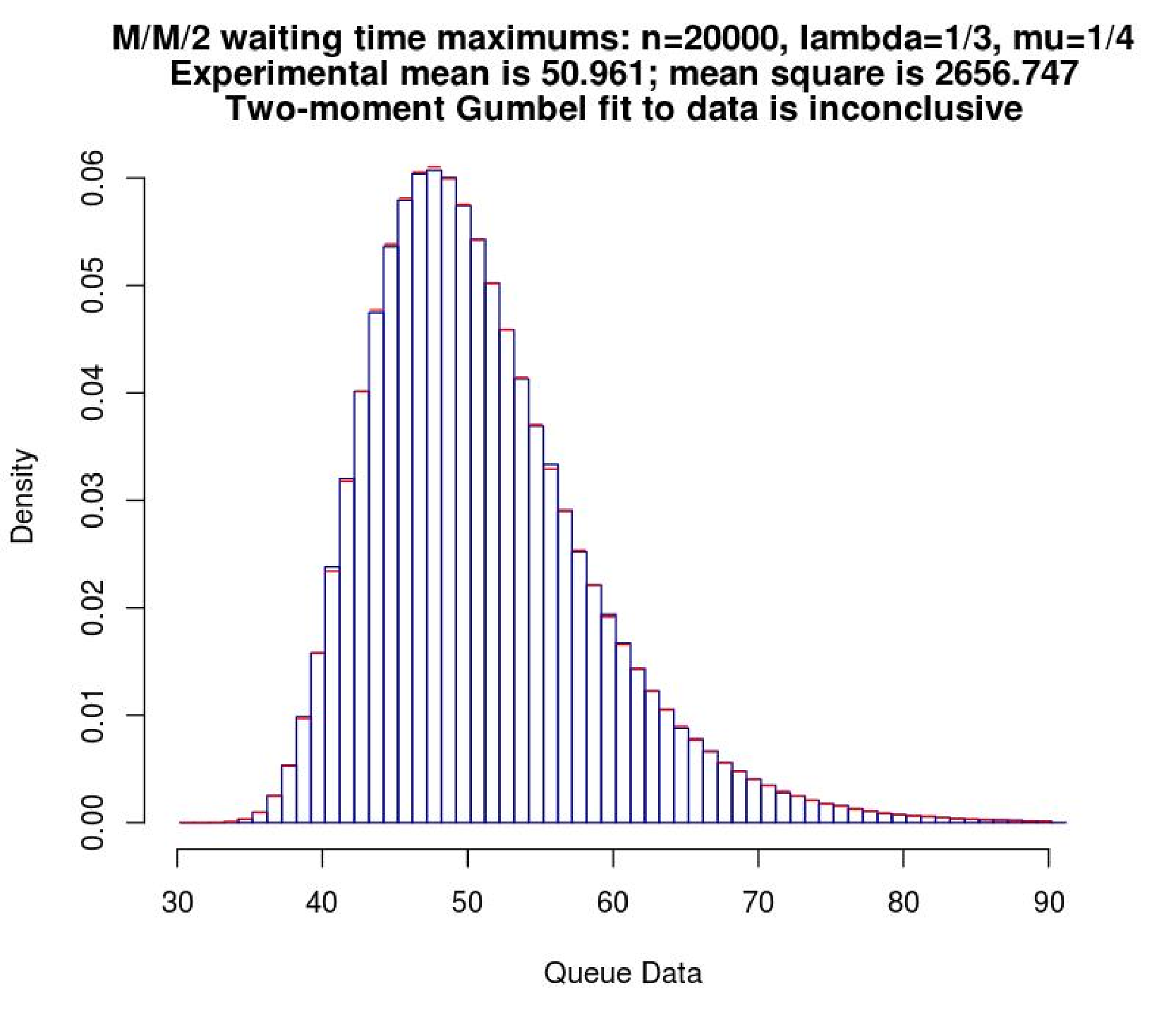}%
\caption{Maximum of $W_{sys}$ over time interval of length $n$, with $2$
servers}%
\end{center}
\end{figure}
\begin{figure}
[ptb]
\begin{center}
\includegraphics[
height=3.0234in,
width=3.4541in
]%
{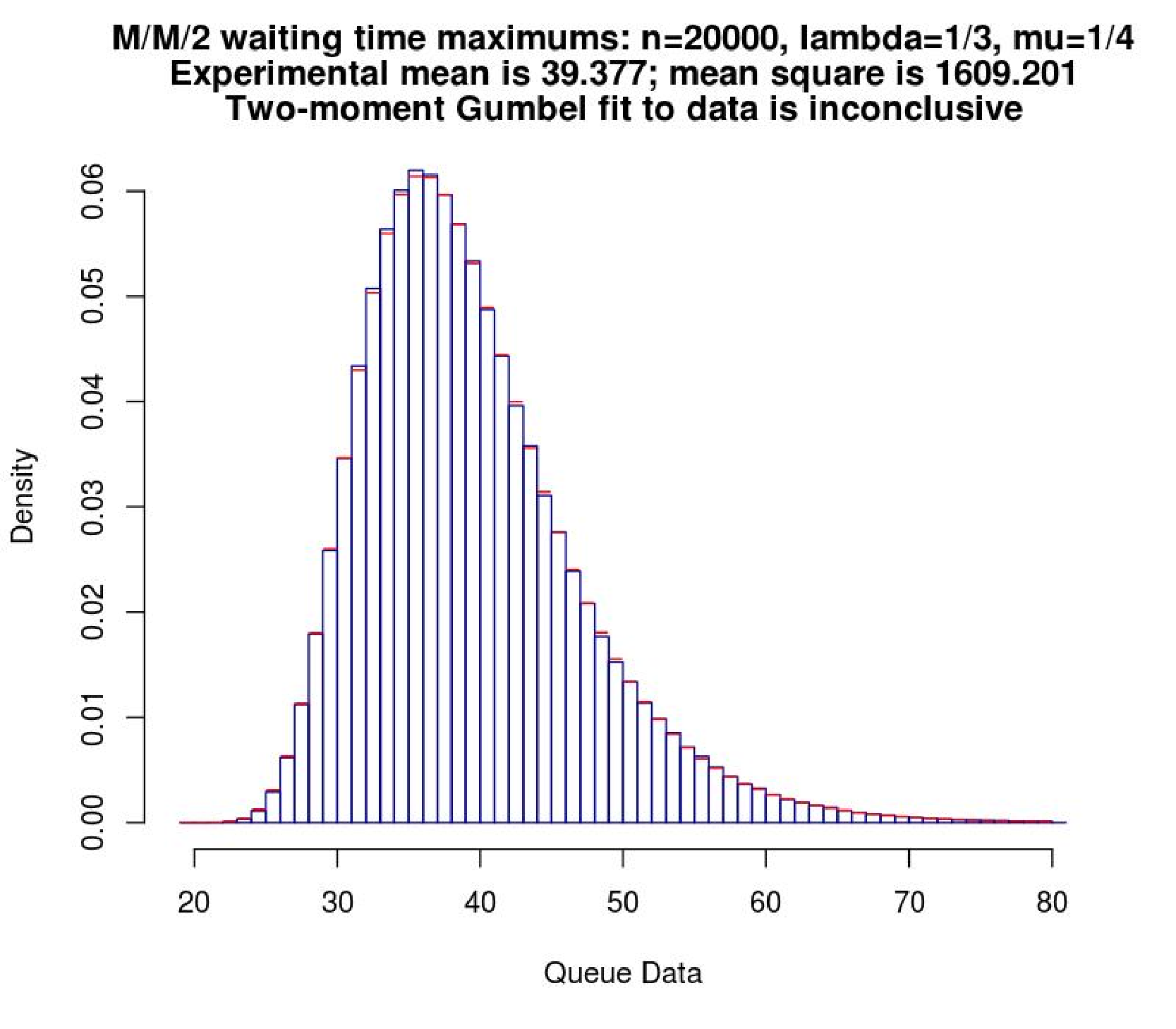}%
\caption{Maximum of $W_{que}$ over time interval of length $n$, with $2$
servers}%
\end{center}
\end{figure}
\begin{figure}
[ptb]
\begin{center}
\includegraphics[
height=3.0234in,
width=3.4541in
]%
{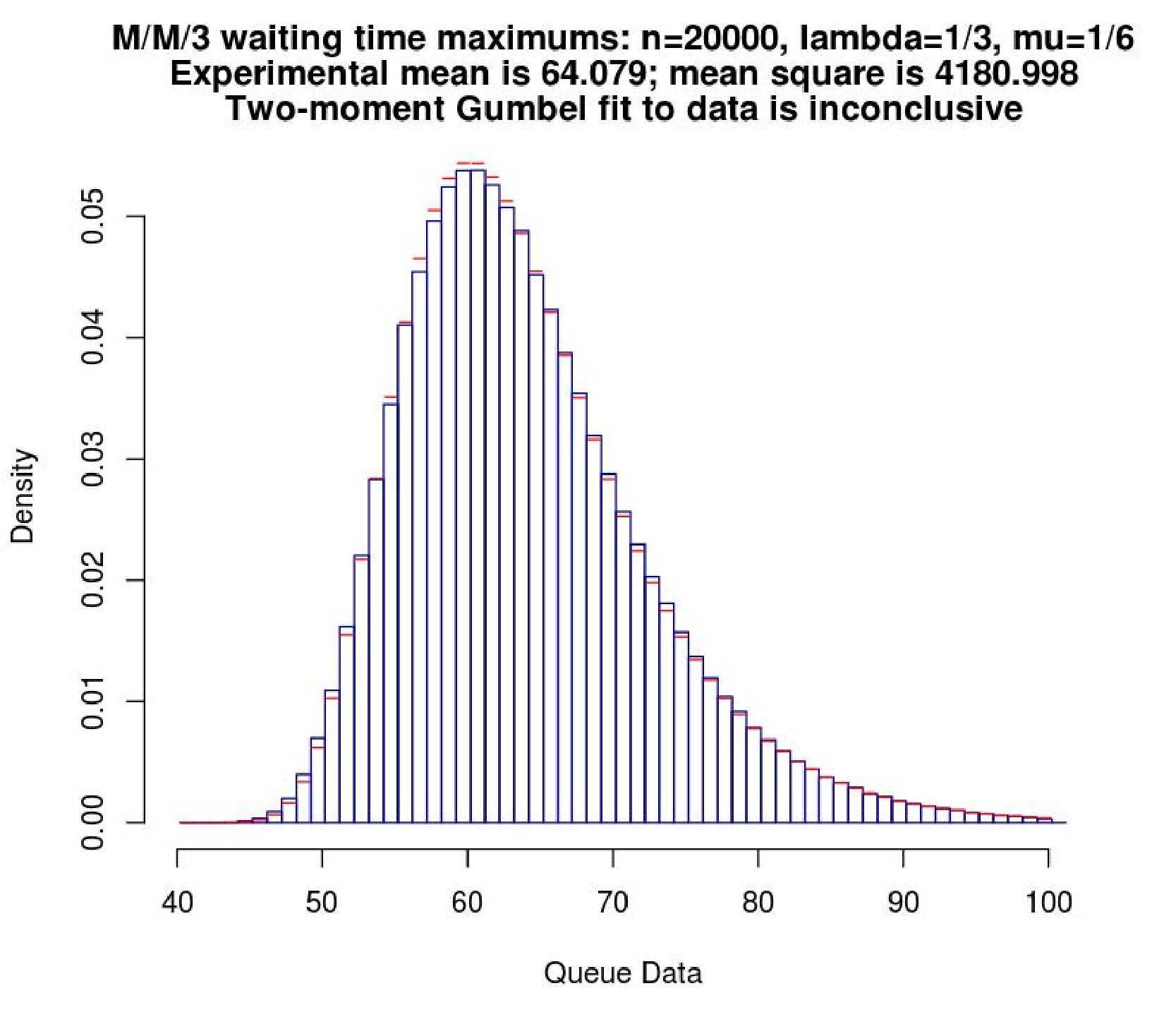}%
\caption{Maximum of $W_{sys}$ over time interval of length $n$, with $3$
servers}%
\end{center}
\end{figure}
\begin{figure}
[ptb]
\begin{center}
\includegraphics[
height=3.0234in,
width=3.4541in
]%
{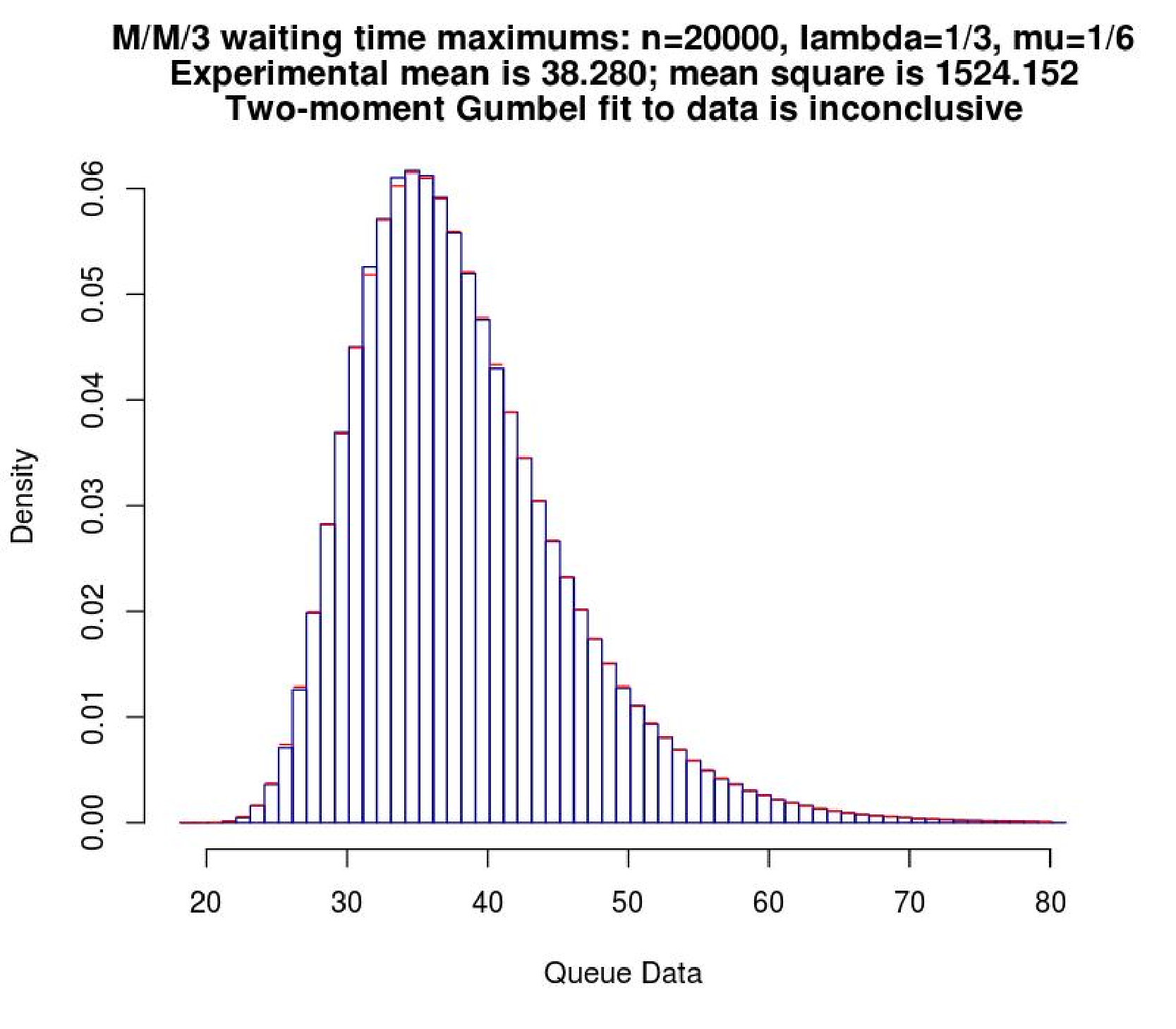}%
\caption{Maximum of $W_{que}$ over time interval of length $n$, with $3$
servers}%
\end{center}
\end{figure}

\pagebreak


\begin{thebibliography}{99}                                                                                               %


\bibitem {Htr-heu}J. J. Hunter, \textit{Mathematical Techniques of Applied
Probability.} Vol. 2: \textit{Discrete Time Models: Techniques and
Applications}, Academic Press, 1983, pp. 189--200; MR0719019.

\bibitem {Fnh-heu}S. Finch, Geo/Geo/$2$ queues and the Poisson clumping
heuristic, arXiv:1902.09272.

\bibitem {Ald-heu}D. Aldous, \textit{Probability Approximations via the
Poisson Clumping Heuristic}, Springer-Verlag, 1989, pp. 1--8, 23--25, 30; MR0969362.

\bibitem {AH-heu}J. R. Artalejo and O. Hern\'{a}ndez-Lerma, Performance
analysis and optimal control of the Geo/Geo/$c$ queue, \textit{Performance
Evaluation} 52 (2003)\ 15--39.

\bibitem {Fi1-heu}S. R. Finch, Euler-Mascheroni constant, \textit{Mathematical
Constants}, Cambridge Univ. Press, 2003, pp. 28--40; MR2003519.

\bibitem {Fi2-heu}S. Finch and G. Louchard, Conjectures about traffic light
queues, arXiv:1810.03906.

\bibitem {Ig1-heu}D. L. Iglehart, Extreme values in the GI/G/$1$ queue,
\textit{Annals Math. Statist.} 43 (1972) 627--635; MR0305498

\bibitem {Ig2-heu}D. L. Iglehart, \textit{Regenerative Simulation for Extreme
Values}, Stanford Univ. Technical Report 43 (1977); http://apps.dtic.mil/docs/citations/ADA047945.

\bibitem {BW-heu}A. W. Berger and W. Whitt, Maximum values in queueing
processes, \textit{Probab. Engrg. Inform. Sci.} 9 (1995) 375--409; MR1365267.

\bibitem {Flt-heu}L. Flatto, The waiting time distribution for the random
order service $M/M/1$ queue, \textit{Annals Appl. Probab.} 7 (1997) 382--409; MR1442319

\bibitem {Lbv-heu}A. V. Lebedev, Maxima of waiting times for the random order
service $M$%
$\vert$%
$M$%
$\vert$%
$1$ queue (in Russian); \textit{Problemy Peredachi Informatsii} v. 41 (2005)
n. 3, 123--127; Engl. transl. in \textit{Problems of Information Transmission}
v. 41 (2005) n. 3, 296--299; MR2163855.

\bibitem {Gm-heu}M. I. Gomes, On maxima of waiting times, \textit{Portugal.
Math.} 39 (1980) 331--339 (1985); MR0776245.

\bibitem {AA-heu}Y. H. Abdelkader and M. Al-Wohaibi, Computing the performance
measures in queueing models via the method of order statistics, \textit{J.
Appl. Math.} 2011, art. 790253; MR2854962.

\bibitem {Alf-heu}A. S. Alfa, \textit{Applied Discrete-Time Queues},
2$^{\text{nd}}$ ed., Springer-Verlag, 2016, pp. 143--150, 227--236; MR3468938.

\bibitem {Fi3-heu}S. Finch, The maximum of an asymmetric simple random walk
with reflection, arXiv:1808.01830.

\bibitem {Fi0-heu}S. Finch and G. Louchard, Traffic light queues and the
Poisson clumping heuristic, arXiv:1810.12058.

\bibitem {Fic-heu}S. Finch, M/M/$c$ queues and the Poisson clumping heuristic, arXiv:1904.04054.

\bibitem {Fi4-heu}S. Finch, M/G/$1$-FIFO queue with uniform service times, arXiv:2206.11108.

\bibitem {Fi5-heu}S. Finch, M/D/$1$ queues with LIFO and SIRO policies, arXiv:2208.09980.%

\begin{tabular}
[c]{lll}
& Steven Finch & \\
& MIT Sloan School of Management & \\
& Cambridge, MA, USA & \\
& \textit{steven\_finch@harvard.edu} &
\end{tabular}

\end{thebibliography}
\end{document}